\documentstyle[11pt,leqno]{article}

\newtheorem{defi}{Definition}[section]
\newtheorem{cor}[defi]{Corollary}

\newtheorem{thm}[defi]{Theorem}
\newtheorem{lem}[defi]{Lemma}

\newtheorem{exa}[defi]{Example}

\newcommand{\D}{\mbox{\tt\normalsize D}}

\renewcommand{\l}{\langle}
\newcommand{\ra}{\rangle}

\newcommand{\be}{\begin{equation}}
\newcommand{\ee}{\end{equation}}

\newcommand{\pf}{ Proof: \ }

\def\eop{\hfill\qquad\rule[-1mm]{2.5mm}{2.5mm}}

\newcommand{\R}{\rm I\kern-.19emR}

\begin{document}

\bibliographystyle{plain}

\title{INVERSE TRIDIAGONAL Z--MATRICES}
\author{
J. J. McDonald
\thanks{Work supported by an NSERC Research Grant.}\\
Department of Mathematics and Statistics\\ University of Regina\\
Regina, Saskatchewan S4S 0A2
\and
R. Nabben
\thanks{Research supported partially by the
Deutsche Forschungsgemeinschaft and NSF Grant DMS-9424346.}\\
Fakult\"at f\"ur Mathematik\\
Universit\"at Bielefeld\\
Postfach 10 01 31 \\
33501 Bielefeld, Germany
\and
M. Neumann
\thanks{Work supported by NSF Grant DMS-9306357.}\\
Department of Mathematics\\
University of Connecticut\\
Storrs, Connecticut 06269--3009
\and
H. Schneider
\thanks{Work supported by NSF Grant DMS-9424346.}\\
Department of Mathematics\\
University of Wisconsin\\
Madison, Wisconsin 53706
\and
M. J. Tsatsomeros$\ ^{*}$\\
Department of Mathematics and Statistics\\ University of Regina\\
Regina, Saskatchewan S4S 0A2}

\date{ \  \\
       \  \\
       \  \\
Dedicated to Robert C. Thompson\\
in memory of his great contributions to linear algebra}

\maketitle
\thispagestyle{empty}

\newpage
\begin{center}
{\bf Abstract}
\end{center}
In this paper, we consider matrices whose inverses are
tridiagonal Z--matrices. Based
on a characterization of symmetric tridiagonal matrices
by Gantmacher and Krein, we
show that a matrix is the inverse of a tridiagonal
Z--matrix if and only if, up to a
positive scaling of the rows, it is the Hadamard product
of a so called weak type $\D$
matrix and a flipped weak type $\D$ matrix whose
parameters satisfy certain quadratic
conditions. We predict from these parameters to which
class of Z--matrices the inverse belongs to. In particular, we give a
characterization of inverse tridiagonal M--matrices. Moreover, we characterize
inverses of tridiagonal M--matrices that satisfy certain row sum criteria.
This leads to the cyclopses that are matrices
constructed from type $\D$ and flipped type $\D$
matrices. We establish some properties of the cyclopses and
provide explicit formulae for the entries of the
inverse of a nonsingular cyclops. We also show that the cyclopses
are the only generalized ultrametric matrices whose inverses are tridiagonal.

\newpage

%********************************************************************
\section{Introduction}
%********************************************************************

In many mathematical problems,
Z--matrices and M--matrices play an important role. It is often useful to
know the properties of their inverses, in particular when the Z--matrices
and the M--matrices have a special combinatorial structure.
In this paper, we investigate the properties
of inverse tridiagonal Z--matrices and M--matrices, i.e., matrices
whose inverses are tridiagonal Z--matrices or M--matrices.
We also highlight some connections between
weak type $\D$ matrices (a class that generalizes type $\D$
matrices as defined by Markham \cite{M}) and inverse tridiagonal
Z--matrices.

First, under the assumption of irreducibility,
we show that a matrix is the inverse of a tridiagonal Z--matrix
if and only if, up to a positive scaling of the rows, it is the
Hadamard product of a weak type $\D$
matrix and a flipped weak type $\D$ matrix whose parameters satisfy certain
quadratic conditions (Theorem \ref{invz}).
This characterization parallels (and is
based on) the characterization of (symmetric) Green matrices
by Gantmacher and Krein \cite{GK}. Further, recalling the classification of
Z--matrices by Fiedler and Markham \cite{FM}, we predict the class $L_s$ of a
tridiagonal Z--matrix based on the parameters of the associated weak
type $\D$ matrices (Theorem \ref{zclass}).
In particular, we find conditions on the parameters so that the inverse is a
tridiagonal M--matrix (Corollary \ref{mgreen}).

Next, we associate type $\D$ matrices with tridiagonal
Z--matrices via the so called cyclopses. These are matrices
that admit a block partition comprising two diagonal blocks that
are of flipped type
$\D$ and of type $\D$, respectively, and two off-diagonal blocks
that have constant entries. We find conditions on the parameters of
the associated type $\D$ matrices and the constant off--diagonal
entries so that the inverse of a cyclops
exists and is a tridiagonal Z--matrix; its nonzero entries are
also found explicitly in terms of the parameters (Theorem \ref{invcyc}).
When a cyclops is a priori nonsingular, we provide necessary and
sufficient conditions so
that its inverse is a tridiagonal Z--matrix (Corollary \ref{zm});
as before we can
predict the class $L_s$ of the tridiagonal Z--matrix (Theorem \ref{czclass}).

Cyclopses (with nonnegative entries) were encountered by the authors
as a special case of the generalized ultrametric matrices (see \cite{NV} and
\cite{MNST1}), which is a class of inverse (row and column diagonally dominant)
M--matrices. We conclude by finding necessary and sufficient conditions so
that
a cyclops is the inverse of a (row and column) diagonally dominant tridiagonal
M--matrix (or equivalently a totally nonnegative generalized ultrametric
matrix)
(see Theorems \ref{igum}, \ref{rgum}, and \ref{ggum}). These results amount to
a characterization of the generalized ultrametric matrices whose inverses are
tridiagonal.

We continue with the precise definitions of the terms mentioned above and
the notational conventions.

%********************************************************************
\section{Preliminaries}
%******************************************************************
We let $e$ denote the all ones vector and $e_j$ the $j$--th standard basis
vector in $\R^n$. Given a positive integer $n$ we let $\l n\ra
=\{1,2,\ldots,n\}$. Let
$\circ$ denote the Hadamard (i.e., entrywise) product of matrices.
For $A = [a_{ij}] \in \R^{n,n}$, by $A(i|j)$
we denote the submatrix of $A$ obtained by deleting the
$i$--th row and the $j$--th column. Given $R,S\subseteq\l n\ra$
we write $A_{RS}$ for
the submatrix of $A$ whose rows and columns are indexed by $R$ and $S$,
respectively.
If $S=\l n\ra\setminus R$ and if $A_{RR}$ is nonsingular,
then the {\em Schur complement} of $A_{RR}$ in $A$ is defined and denoted by
\[ A/A_{RR}=A_{SS}-A_{SR}(A_{RR})^{-1}A_{RS}. \]
It is well known that $\mbox{det}A=\mbox{det}A_{RR}\mbox{det}(A/A_{RR})$.

We call $A=[a_{ij}]\in\R^{n,n}$ a Z--matrix if $a_{ij}\leq 0$
for all $i\neq j$. For any nonnegative integer $s\leq n$
we denote by $L_s$ the set of all matrices $A=tI-B\in\R^{n,n}$,
where $B$ is an entrywise nonnegative matrix and where $\rho_s(B)\leq t <
\rho_{s+1}(B)$. Here $\rho_s(B)$ denotes the maximum among the
spectral radii of all
the $s\times s$ principal submatrices of $B$ (we take $\rho_0=-\infty$ and
$\rho_{n+1}=\infty)$. In particular, $A$ is an {\em M--matrix} if it can be
written as $A=tI-B$, where $B$ is an entrywise nonnegative matrix and
$\rho(B):=\rho_n(B)\leq t$.

The next theorem, found in \cite{N} and \cite{Sm}, is a
characterization of the nonsingular Z--matrices in $L_s$.

\begin{thm}
\label{lsz}
Let $A \in \R^{n,n}$ be nonsingular Z--matrix. Then $A\in L_s$ if
and only if one of the following alternative cases a) or b) holds:
\begin{description}
\item[{\bf a)}]
\begin{description}
\item[{\bf \hspace*{.1in} (i)}]
\mbox{\em det}$A<0$,
\item[{\bf (ii)}]
all principal minors of $A^{-1}$ of order greater
than or equal to $n-s$ are nonpositive, and
\item[{\bf (iii)}]
there exists a positive principal minor of $A^{-1}$ of order $n-s-1$.
\end{description}
\item[{\bf b)}]
\begin{description}
\item[{\bf\hspace*{.1in} (i)}]
\mbox{\em det}$A>0$,
\item[{\bf (ii)}]
all principal minors of $A^{-1}$ of order greater
than or equal to $n-s$ are nonnegative, and
\item[{\bf (iii)}]
there exists a negative principal minor of $A^{-1}$ of order $n-s-1$.
\end{description}
\end{description}
\end{thm}

Markham defined in \cite{M} type $\D$ matrices as follows:
$A = [a_{ij}] \in \R^{n,n}$ is of {\em type $\D$}
(or {\em a type $\D$ matrix}) if
\begin{eqnarray*} \label{deftypd}
a_{ij} = \left\{
\begin{array}{c}
a_i, \ \ i \leq j, \\
a_j, \ \ i > j,
\end{array} \right.
\quad \mbox{where} \ a_n > a_{n-1} > \ldots > a_1. \end{eqnarray*}
We refer to the $a_i$ ($i=1,2,\ldots,n$) as the {\em parameters} of
$A$. We also consider similarly constructed matrices, without constraints
on the parameters $a_{i}$, to which we refer as of {\em weak type $\D$}.
Moreover, we call $A$ a {\em flipped type $\D$ matrix}
(resp., a {\em flipped weak type $\D$ matrix})
if $PAP^T$ is a type $\D$ matrix (resp., a weak type $\D$
matrix), where $P$ is the permutation that reverses the order of the indices
$1,2,\ldots,n$. We enumerate the parameters of a weak type $\D$ matrix,
as well as the parameters of a flipped weak type $\D$ matrix in a way
such that the $i$--th parameter is equal to the $i$--th diagonal
entry of the matrix. To illustrate these definitions and the
relevant notation, let
\[ A=\left[\matrix{-1&-1&-1\cr
                   -1&\ 2&\ 2\cr
                   -1&\ 2&\ 3\cr}\right], \ \
B=\left[\matrix{-3&2&1\cr
               \ 2&2&1\cr
               \ 1&1&1\cr}\right].
\]
Then $A$ is of type $\D$ with parameters $a_i$ given by
$(-1,2,3)$ and $B$ is of flipped weak type $\D$ with parameters $b_i$
given by $(-3,2,1)$.

Gantmacher and Krein defined in \cite{GK} a {\em Green matrix}
to be a matrix $G \in \R^{n,n}$ such that $G = A \circ B$, where $A$
is a weak type $\D$ matrix, $B$ is a flipped weak type $\D$ matrix.
The name Green matrix is
not the only name for these matrices. Originally Gantmacher and Krein
called such matrices {\em einpaarig} or {\em matrix of a couple}.
Moreover, Markham defined the
type $\D$ matrices as a special case of Green matrices.

In our discussion, we shall also refer to the following matrices that were
introduced in \cite{NV} and \cite{MNST1}. We say $C=[c_{ij}]\in \R^{n,n}$
is a {\em generalized ultrametric matrix} if
\begin{description}
\item[(i)] $C$ is entrywise nonnegative,
\item[(ii)] $c_{ii} \geq \max \{c_{ij},c_{ji} \}$ for all $i,j \in \l n \ra$,
\item[(iii)] every subset of
$\l n \ra$ with three distinct elements has a labeling $\{i,j,k\}$ such that
\begin{description}
\item[(a)] $c_{ij}=c_{ik},$
\item[(b)] $c_{ji}=c_{ki},$
\item[(c)] $\min \{ c_{jk},c_{kj}\} \geq \min \{c_{ji},c_{ij}\},
$ \item[(d)] $\max \{ c_{jk},c_{kj}\} \geq \max \{c_{ji},c_{ij}\}.$
\end{description}
\end{description}
In the aforementioned papers, it is shown that if a generalized
ultrametric matrix is nonsingular then its inverse is a row and column
diagonally dominant M--matrix.

Next, we introduce a class of matrices
constructed from type $\D$ matrices; as we show in
Section 4, it contains matrices that are under certain additional conditions
are inverse tridiagonal Z--matrices. Let $C\in\R^{n,n}$ and let $m\leq n$ be
a nonnegative integer. We call $C$ a {\em cyclops with eye $m+$} if
\begin{eqnarray}
\label{cyclops}
C = \left[ \begin{array}{cc}
C_{11} & b_1E_{12} \\
b_2E_{21} & C_{22}
\end{array} \right],
\end{eqnarray}
where $C_{11}$ is a $m\times m$ flipped type $\D$ matrix and $C_{22}$ is a
$(n-m)\times (n-m)$ type $\D$ matrix, viz., \begin{eqnarray*}
C_{11} = \left[ \begin{array}{ccccc}
a_1 & a_2 & \ldots & a_{m-1} & a_m \\
a_2 & a_2 & \ldots & a_{m-1} & a_m \\
\vdots & \vdots & & \vdots & \vdots \\
\vdots & \vdots & \ddots & \vdots & \vdots \\
a_{m-1} & a_{m-1} & \ldots & a_{m-1} & a_m \\ a_{m} & a_m & \ldots & a_{m}
& a_m
\end{array} \right],
\quad
C_{22} = \left[ \begin{array}{ccccc}
a_{m+1} & a_{m+1} & \ldots & a_{m+1} & a_{m+1} \\
a_{m+1} & a_{m+2} & \ldots & a_{m+2} & a_{m+2} \\
\vdots & \vdots & & \vdots & \vdots
\\
\vdots & \vdots & \ddots & \vdots & \vdots \\
a_{m+1} & a_{m+2} & \ldots & a_{n-1} & a_{n-1} \\ a_{m+1} & a_{m+2} & \ldots &
a_{n-1} & a_n \\ \end{array} \right]
\end{eqnarray*}
with
\begin{eqnarray}
\label{inequa}
a_1 > a_2 > \ldots > a_m \quad \mbox{and} \quad a_{n}
>a_{n-1} > \ldots > a_{m+1},
\end{eqnarray}
and where $E_{12}$ and $E_{21}$ are all ones matrices of appropriate sizes.
We refer to the $a_i$ ($i=1,2,\ldots,n$) and $b_1,b_2$ as the
{\em parameters} of the cyclops $C$.

In the remainder of this paper, when we refer to a type $\D$ matrix,
a weak type $\D$ matrix, a Green matrix, or a cyclops, we assume that
the reader recalls the notation and the associated parameters indicated in
this section.

%********************************************************************
\section{Hadamard Products of weak type D Matrices}
%*****************************************************************
Gantmacher and Krein proved the following results.

\begin{thm}
{\em (Gantmacher and Krein \cite{GK})}
\label{gantkrein}
Let $G \in \R^{n,n}$ be symmetric. Then the following are equivalent:
\begin{description}
\item[{\bf (i)}] $G$ is a nonsingular Green matrix.
\item[{\bf (ii)}] $G^{-1}$ is an irreducible tridiagonal matrix.
\end{description}
\end{thm}

\begin{lem} \label{lemgk}
{\em (Gantmacher and Krein \cite{GK})}
\label{lemmadet}
Let $G \in \R^{n,n}$ be a Green matrix with associated parameters
$a_i, b_i$. Let $h_i := a_ib_{i-1} - a_{i-1}b_i$ for $i=2,3,\ldots,n$. Then
\begin{eqnarray*}
\mbox{\em det}G=a_1 b_n \prod_{i=2}^{n} h_i. \end{eqnarray*}
Moreover,
\begin{eqnarray*}
\mbox{\em det} G(i|j) = \left\{
\begin{array}{rl}
\mbox{\em det}G/h_{i+1} & \mbox{if} \ |i-j| = 1 \\ 0 & \mbox{if } \ |i-j| > 1.
\end{array} \right.
\end{eqnarray*}
\end{lem}

We proceed by characterizing inverse
tridiagonal Z--matrices in the spirit of Theorem \ref{gantkrein}.

\begin{thm}
\label{invz}
Let $C \in \R^{n,n}$ be nonsingular and irreducible.
Then the following are equivalent:
\begin{description}
\item[{\bf (i)}] $C^{-1}$ is a tridiagonal Z--matrix.
\item[{\bf (ii)}] There exists a positive diagonal matrix $D\in\R^{n,n}$
such that $DC=A\circ B$, where $A$ is of weak type $\D$ with parameters
$a_i$, and $B$ is of flipped weak type $\D$ with parameters $b_i$,
such that $a_ib_{i-1}-a_{i-1}b_i>0$ for all $i=2,3,\ldots,n$.
\end{description}
\end{thm}

\pf
Let $h_i:= a_ib_{i-1} - a_{i-1}b_i$ for $i=2,3,\ldots,n$ and suppose that (i)
holds.  As is well known (see \cite{ES} and \cite{PY}), there exists a
positive
diagonal matrix $D^{-1}$ such that $C^{-1}D^{-1}$ is symmetric.
Thus, by Theorem \ref{gantkrein},
$DC$ is a Green matrix. We also have that $C^{-1}D^{-1} = [\gamma_{ij}]$,
where
\[
\gamma_{ij} = (-1)^{i+j}\mbox{det}((DC)(j|i))/\mbox{det}(DC). \]
Hence, by Lemma \ref{lemmadet}, the superdiagonal entries are \be
\label{gamma}
\gamma_{i,i+1} = \frac{-1}{h_{i+1}} \ \ (i=1,2,\ldots,n). \ee
Thus all $h_i$ are positive. Conversely, if (ii) holds, (\ref{gamma}) and
Theorem \ref{gantkrein} imply that $C^{-1}D^{-1}$ is a tridiagonal Z--matrix.
Hence $C^{-1}$ is a tridiagonal Z--matrix. \eop

In the following theorem we determine the class $L_s$ to which an inverse
tridiagonal Z--matrix belongs.

\begin{thm}
\label{zclass}
Let $C \in \R^{n,n}$ be an irreducible inverse tridiagonal
Z--matrix. Let $D$, $A$, and $B$ be as in condition (ii)
of Theorem \ref{invz}. Then the following hold:
\begin{description}
\item[{\bf (i)}] If $\mbox{\em det}C < 0$ then $C^{-1}
\in L_s$ with $s = \min\{t-2,n-r-2\}$;
\item[{\bf (ii)}] if $\mbox{\em det}C > 0$ then $C^{-1}
\in L_s$ with $s = \min\{t-2,n-q-2\}$,
\end{description}
where
\begin{eqnarray*}
t & = & \left\{\begin{array}{l}
n+2 \ \ \mbox{if $a_ib_{i+j} - a_{i+j}b_i > 0 $ for all $i,j$}\\
\min\{j \in \l n\ra | \
        \mbox{there exists an} \ i \in \l n\ra  \
\mbox{with} \ a_ib_{i+j} - a_{i+j}b_i< 0\} \ \mbox{otherwise,}
\end{array}
\right. \\
 r & = & \left\{\begin{array}{l}
-1 \ \ \mbox{if $a_ib_{i+j}\leq 0$  for all $i,j$}\\
\max\{j \in \l n\ra  | \ \mbox{there exists an} \ i \in \l n\ra  \
\mbox{with} \ a_ib_{i+j}>0 \} \ \ \mbox{otherwise, \ }
\end{array}
\right. \\
q & = & \left\{\begin{array}{l}
-2,\ \ \mbox{if $a_ib_{i+j}>0$  for all $i,j$}\\
\max\{j \in \l n\ra  | \ \mbox{there exists an} \ i \in \l n\ra  \
\mbox{with} \ a_ib_{i+j}< 0\}, \ \ \mbox{otherwise.}
\end{array}
\right.
\end{eqnarray*}
\end{thm}

\pf
For $i,j \in \l n\ra$ with $i>j$, define $h_{ij}=a_ib_j - a_jb_i.$
Since $C$ is an inverse tridiagonal Z--matrix we have,
by Theorem \ref{invz}, that $h_{i,i-1} > 0$ for all $i=2,3,\ldots,n$.
Moreover, by Lemma \ref{lemmadet}, we have that
\begin{eqnarray}
\label{dethier}
\mbox{det}(DC) = a_1 b_n \prod_{i=2}^{n} h_{i,i-1},
\end{eqnarray}
where $D$ is the positive diagonal matrix chosen in the proof of
Theorem \ref{invz}. We will proceed by considering the signs of the
principal minors of $DC$ and by applying Theorem \ref{lsz}.
Since principal submatrices of Green matrices are also Green matrices,
the principal minors of $DC$ are given by formulae similar to (\ref{dethier}),
and their signs are determined by the corresponding quantities $h_{ij}$
and $a_i$ and $b_i$.

First, suppose $\mbox{det} C<0$, i.e., $a_1b_n<0$.
Without loss of generality, we can assume that $a_1 > 0$ and $b_n < 0$.
When $t \neq n + 2$, the definition of $t$ and (\ref{dethier})
imply that there exists a principal submatrix of order $n - (t - 1)$ of
$DC$ with positive determinant. (This principal submatrix is obtained by
deleting rows and columns $i+1,\ldots,i+j-1$, where $i,j$ are the minimal
indices in the definition of $t$.)
For all principal submatrices of order greater than $n - (t - 1)$, the
relevant $h_{i,j}$ appearing in the determinantal
formula of Lemma \ref{lemmadet} are positive. It is also clear
from the definition of $r$ that there exists an $(r+1)\times(r+1)$
principal submatrix with positive determinant. Moreover, all principal
submatrices of order $\tilde n$ with $ \tilde n> r+1$ satisfy
$\tilde a_1\tilde b_{\tilde n} \leq 0,$
where $\tilde a_1$ and $\tilde b_{\tilde n}$ are the corresponding parameters.
>From these cases, we obtain the following: if $s = \min\{t-2,n-r-2\}$,
then there exists a principal minor of order $n-s-1$ that is positive.
Also, all principal minors of order greater than $n-s-1$ are nonpositive.
Thus, by Theorem \ref{lsz}, $C^{-1}\in L_s$, showing (i).
Similarly we obtain (ii).
\eop

It is shown in \cite{N} that if $C^{-1}\in L_s$ and
$\lfloor\frac{n}{2}\rfloor\leq s<n$, then det$C < 0$ .
For inverse tridiagonal Z--matrices, this result can be established by
considering the changes of the signs in the sequences of the parameter
$a_i$ and $b_i$. With $q$ as in Theorem \ref{zclass}, if det$C>0$ and $C$
is not entrywise nonnegative, one obtains that
$q+1\geq\lfloor\frac{n}{2}\rfloor$.

\begin{exa}
{\em In the following examples we apply Theorem \ref{zclass}.\\
{\bf (i)} Consider the matrix
\begin{eqnarray*}
C = \left[ \begin{array}{rrrr}
-2 & -2 & -2 &-2 \\
-2 & -1 & -1 & -1 \\
-2 & -1 & 2 & 2\\
-2 & -1 & 2 & 3
\end{array} \right]
\end{eqnarray*}
for which we can write $DC=A\circ B$, where $D=I$, $A$ is of
type D with parameters $a_i$ given by $(-2,-1,2,3)$, and $B$ is of
(flipped) weak type D with parameters $b_i$ given by $(1,1,1,1)$.
Notice that $a_ib_{i-1}-a_{i-1}b_i>0$ for all $i=2,3,4$.
Moreover, $\mbox{det}C<0$, $t = 6$ and $r = 1$. Thus $C^{-1}\in L_1$.\\
{\bf (ii)} Let
\begin{eqnarray*}
C = \left[ \begin{array}{rrrr}
-24 & -20 & -16 & -4\\
-20 & -10 & -8 & -2 \\
-16 & -8 & 16 & 4 \\
-4 & -2 & 4 & 5
\end{array} \right] =
\left[ \begin{array}{rrrr}
-4 & -4 & -4 & - 4\\
-4 & -2 & -2 & -2 \\
-4 & -2 & 4 & 4\\
-4 & -2 & 4 & 5
\end{array} \right] \circ
\left[ \begin{array}{rrrr}
6 & 5 & 4 & 1 \\
5 & 5 & 4 & 1\\
4 & 4& 4 & 1\\
1 & 1& 1& 1
\end{array} \right].
\end{eqnarray*}
$C$ is the Hadamard product of a weak type D and a flipped weak type D
matrix with parameters $(-4,-2,4,5)$ and $(6,5,4,1)$, respectively.
By Theorem \ref{invz}, $C^{-1}$ is a tridiagonal Z--matrix
and $\mbox{det}C<0$. Moreover, $r=1$ and $t=2$ and hence $C^{-1}\in L_0$.
}
\end{exa}

The results above yield the following characterization of inverse tridiagonal
M--matrices.

\begin{cor}
\label{mgreen}
Let $C \in \R^{n,n}$ be nonsingular and irreducible.
Then the following are equivalent:
\begin{description}
\item[{\bf (i)}] $C^{-1}$ is a tridiagonal M--matrix
\item[{\bf (ii)}]
There exists a positive diagonal matrix $D \in \R^{n,n}$ such that
$DC=A\circ B$, where $A$ is of a weak type $\D$ with parameters $a_i$,
and $B$ is of flipped weak type $\D$ with parameters $b_i$, such that all
the parameters have the same sign and
\begin{eqnarray*} 0 < \frac{a_1}{b_1} <
\frac{a_2}{b_2} < \ldots < \frac{a_n}{b_n}. \end{eqnarray*}
\end{description}
\end{cor}

\pf

{\bf (i) implies (ii):}
As $C^{-1}\in L_n$, we have that det$C>0$ and that
$C^{-1}$ is a tridiagonal Z--matrix. The implication now follows
from Theorems \ref{invz} and \ref{zclass}.

{\bf (ii) implies (i):} \
By Theorem \ref{invz}, $C^{-1}$ is a tridiagonal Z--matrix. Since $C$ is
entrywise positive and since every inverse positive
Z--matrix is an M--matrix, (i) holds.  \eop

%*************************************************************
\section{Cyclopses}
%*************************************************************
In this section, we consider inverse tridiagonal Z--matrices that satisfy
certain row sum and column sum criteria. This leads to a new class of
matrices that we have defined as cyclopses in Section 2. 
We begin with some auxiliary results.
\begin{lem}
{\em (\cite[Observation 3.8]{N})}
\label{rsumd}
Let $A\in\R^{n,n}$ be of type $D$ with parameters $a_i$ and such that
$a_1\neq 0$.
Then
$$A^{-1}e=\left[\frac{1}{a_1},0,0,\ldots,0\right]^T.$$
\end{lem}
>From Lemma \ref{lemgk}, for a type $\D$ matrix with parameters $a_i$,
we obtain that
\begin{eqnarray} \label{detd}
\mbox{\mbox{det}}C=a_1  \prod_{j=2}^n(a_j-a_{j-1}).
\end{eqnarray}

\begin{lem}
\label{schurcyc}
Let $C\in\R^{n,n}$ be a cyclops with eye $m+$ parameters $a_i, b_1, b_2$.
Suppose that $a_m\neq 0$ and $a_{m+1}\neq 0$. Then
$$
C/C_{11}=C_{22}-\frac{b_1b_2}{a_m}E_{1}, \ \ \
C/C_{22}=C_{11}-\frac{b_1b_2}{a_{m+1}}E_{2},
$$
where $E_1, \ E_2$ are all ones matrices of appropriate sizes.
\end{lem}

\pf
Follows from Lemma \ref{rsumd} and the fact that
$C/C_{11}=C_{22}-b_1b_2E_{21}C_{11}^{-1}E_{12}$. \eop

\begin{thm}
\label{detcyc}
Let $C\in\R^{n,n}$ be a cyclops
with eye $m+$ and parameters $a_i, b_1, b_2$. Then
\begin{eqnarray*}
\mbox{\mbox{\em det}}C=
(a_m a_{m+1}-b_1 b_2)\prod_{j=1}^{m-1}(a_{m-j}-a_{m-j+1})
\prod_{j=m+2}^n(a_j-a_{j-1}).
\end{eqnarray*}
\end{thm}

\pf
\begin{description}
\item[Case I ($a_m\neq 0$ or $a_{m+1}\neq 0$):] if $a_m\neq 0$, then by Lemma
\ref{detd}, $C_{11}$ is nonsingular and
$$\mbox{det}C=\mbox{det}C_{11}\mbox{det}(C/C_{11}) $$
$$	=a_{m}\prod_{j=1}^{m-1}(a_{m-j}-a_{m-j+1}) \
\mbox{det}(C_{22}-b_{1}b_{2}/a_{m}E_1) $$
$$	=a_{m}\prod_{j=1}^{m-1}(a_{m-j}-a_{m-j+1})
(a_{m+1}-\frac{b_1b_2}{a_m})
\prod_{j=m+2}^{n}((a_{j}-\frac{b_1b_2}{a_{m}} )-
(a_{j-1}-\frac{b_1b_2}{a_{m}} )) $$
$$ =(a_m a_{m+1}-b_1 b_2)\prod_{j=1}^{m-1}(a_{m-j}-a_{m-j+1})
\prod_{j=m+2}^n(a_j-a_{j-1}).$$
If $a_{m+1}\neq0$ the result follows in a similar manner.

\item[Case II ($a_m=a_{m+1}=0$):]
if $a_m=a_{m+1}=0$ and either $b_1=0$ or $b_2=0$, then $C$ has a row of zeros
(and thus zero determinant) and the result follows. Assume that $b_1\neq 0$ and
$b_2\neq 0$. Let $R=\{m, m+1\}$, $S=\{1,2,\ldots,m-1\},$
$T=\{m+2,m+3,\ldots,n\}$ and
$U=S\cup T$. Then $C_{RR}$ is nonsingular since $a_{m-1}>a_m=0$. Hence
$$
\mbox{det}C=\mbox{det}C_{RR}\mbox{det}(C/C_{RR}) =-b_1b_2 \
\mbox{det}\left[\begin{array}{cc} C_{SS}& 0 \\ 0&C_{TT}
\end{array}\right]
$$
$$
=-b_1b_2(a_{m-1}a_{m+2})\prod_{j=2}^{m-1}(a_{m-j}-a_{m-j+1})
\prod_{j=m+3}^n(a_j-a_{j-1})
$$
$$
=(a_ma_{m+1}-b_1b_2)\prod_{j=1}^{m-1}(a_{m-j}-a_{m-j+1})
\prod_{j=m+2}^n(a_j-a_{j-1}).
$$
\end{description}
\eop

The following is an immediate consequence of the above theorem.
\begin{cor}
\label{singcyc}
Let $C\in\R^{n,n}$ be a cyclops with eye $m+$ and
parameters $a_i, b_1, b_2$. Then $C$ is nonsingular if and only if
$a_ma_{m+1}-b_1b_2\neq 0$.
Moreover, \mbox{\em sgn}$(\det C) =$ \mbox{\em sgn}$(a_ma_{m+1}-b_1b_2)$.
\end{cor}

Next, we shall give explicit formulae for the entries of the inverse of a
nonsingular cyclops.  We first need another result on type $\D$
matrices proved in \cite{N1}. We denote by $\otimes$ the Kronecker
product of matrices.

\begin{thm}
\label{invdex}
Let $A$ be a nonsingular matrix of type $\D$ with parameters
$a_i$. Then the inverse of $A$ is given by
\be
\label{inva}
A^{-1} = \sum_{i=1}^n v^{(i)}(v^{(i)})^T \otimes (a_{i} - a_{i-1})^{-1}, \ee
with $a_0\equiv 0$.
Here the vectors $v^{(i)} = [v^{(i)}_j]$ are defined as \be
\label{vis}
v^{(i)}_j = \left\{\begin{array}{ll}
-1 \ \mbox{if} \ \ j=i-1\\
\ 1 \ \ \mbox{if} \ \ j=i\\
\ 0 \ \ \mbox{otherwise.}
\end{array} \right.
\ee
\end{thm}
It follows that the entries $\alpha_{ij}$ of the inverse of a
type $\D$ matrix $A \in \R^{n,n}$ are zero except for
\[
\alpha_{11} = \frac{1}{a_1} + \frac{1}{a_2-a_{1}}, \ \
\alpha_{nn} = \frac{1}{a_n-a_{n-1}}, \]
\[
\alpha_{ii} = \frac{1}{a_i-a_{i-1}} + \frac{1}{a_{i+1}-a_i} \ \
(i=2,3,\ldots,n-1),
\]
\[
\alpha_{i,i+1} = \alpha_{i+1,i} = - \frac{1}{a_{i+1}-a_i}.
\]
Similarly, for a flipped type $\D$ matrix $A$ we have
\[
\alpha_{11} = \frac{1}{a_{1}-a_2}, \ \ \alpha_{nn} = \frac{1}{a_n} +
\frac{1}{a_{n-1}-a_{n}},
\]
\[
\alpha_{ii} = \frac{1}{a_i-a_{i+1}} + \frac{1}{a_{i-1}-a_i} \ \
(i=2,3,\ldots,n-1),
\]
\[
\alpha_{i,i+1} = \alpha_{i+1,i} = - \frac{1}{a_i-a_{i+1}}. \]

\begin{thm}
\label{invcyc}
Let $C \in \R^{n,n}$ be a cyclops with eye $m+$ and parameters $a_i, b_1,
b_2$.
Suppose that $a_ma_{m+1}-b_1b_2\neq 0$. Then $A=C^{-1}=[\alpha_{ij}]$ exists
and is a tridiagonal matrix with entries given by
\[
\alpha_{11} = \frac{1}{a_{1}-a_2},
\]
\[
\alpha_{ii} = \frac{1}{a_i-a_{i+1}} + \frac{1}{a_{i-1}-a_i} \ \
(i=2,3,\ldots,m-1),
\]
\[
\alpha_{mm} = \frac{a_{m+1}}{a_ma_{m+1}-b_1b_2} + \frac{1}{a_{m-1}-a_{m}},
\]
\[
\alpha_{i,i+1} = \alpha_{i+1,i} = - \frac{1}{a_i-a_{i+1}} \ \
(i=1,2,\ldots,m-1),
\]
\[
\alpha_{m+1,m+1} = \frac{a_m}{a_ma_{m+1}-b_1b_2} +
\frac{1}{a_{m+2}-a_{m+1}}, \]
\[
\alpha_{ii} = \frac{1}{a_i-a_{i-1}} + \frac{1}{a_{i+1}-a_i} \ \
(i=m+2, m+3,\ldots,n-1)
\]
\[
\alpha_{nn} = \frac{1}{a_n-a_{n-1}},
\]
\[
\alpha_{i,i+1} = \alpha_{i+1,i} = - \frac{1}{a_{i+1}-a_i} \ \ (i=m+1,
m+2,\ldots,n-1).
\]
Moreover,
\[
\alpha_{m,m+1} = -\frac{b_1}{a_ma_{m+1}-b_1b_2}, \ \ \alpha_{m+1,m} =
-\frac{b_2}{a_ma_{m+1}-b_1b_2}. \]
\end{thm}

\pf
Recall that our assumption that $a_ma_{m+1}-b_1b_2\neq 0$ is equivalent to $C$
being invertible.
\begin{description}
\item[Case I ($a_m\neq 0$ and $a_{m+1}\neq 0$):] if $a_m\neq 0$ then
$C_{11}$ is a nonsingular flipped type $\D$ matrix and hence it
is the inverse
of a tridiagonal Z--matrix (by the results in \cite{N}).
If $a_{m+1}\neq 0$, then
$C_{22}$ is a nonsingular type $\D$ matrix and hence
it is also the inverse of a tridiagonal Z--matrix.
By Lemma \ref{schurcyc}, $C/C_{11}$ is of type $\D$ and
$C/C_{22}$ is of flipped type $\D$. Moreover, by (\ref{detd}) applied to
$C/C_{11}$ and $C/C_{22}$ and since $a_ma_{m+1}-b_1b_2\neq 0$,
both Schur complements
are nonsingular and thus (using formulas from \cite[(10), p. 773]{BS})
$$
A=\left[\begin{array}{cc} (C/C_{22})^{-1}&
-b_1C_{11}^{-1}E_{12}(C/C_{11})^{-1}\\
-b_2C_{22}^{-1}E_{21}(C/C_{22})^{-1} &(C/C_{11})^{-1} \end{array}\right].
$$
Since $C/C_{22}$ is of flipped type $\D$,
its inverse is a tridiagonal Z--matrix
with all row sums zero except the $m$--th (last).
Since $C/C_{11}$ is of type {\cal
D}, its inverse is a tridiagonal Z--matrix with all
row sums zero except the first,
which corresponds to the $(m+1)$--st row of $A$.

Now one can easily get the entries of
$(C/C_{22})^{-1}$ and $(C/C_{11})^{-1}$ using
Theorem \ref{invdex}. Furthermore, it follows that
$$
-b_1C_{11}^{-1}E_{12}(C/C_{11})^{-1}
=- \frac{b_1}{a_m} e_me_1^T(C/C_{11})^{-1} $$
$$
=- \frac{b_1}{a_m(a_{m+1}-\frac{b_1b_2}{a_m})} e_me_1^T $$
$$
=- \frac{b_1}{a_ma_{m+1}-b_1b_2} e_me_1^T. $$
Similarly,
$$-b_2C_{22}^{-1}E_{21}(C/C_{22})^{-1}
=-\frac{b_2}{a_ma_{m+1}-b_1b_2} e_1e_{n-m}^T. $$
This establishes the result in Case I.

\item[Case II ($a_m=0$ or $a_{m+1}=0$):]
if $a_m=0$ or $a_{m+1}=0$ (or both),
form a new cyclops from $C$ by replacing the $m$--th or the $(m+1)$--st
parameter
(or both) by real numbers that approach zero. The result then follows
from Case I  and continuity. \eop
\end{description}

Notice that if $C$ is as in the previous theorem, then all row sums and column
sums of $A=C^{-1}$ are zero, except at least one of the $m$--th
or the $(m+1)$--st (for otherwise $C$ would be singular).

\begin{cor}
\label{zm}
Let $C\in\R^{n,n}$ be a nonsingular cyclops with eye $m+$ and parameters
$a_i, b_1, b_2$. Then $C^{-1}$ is a Z--matrix if and only if the following two
conditions hold:
\begin{description}
\item[{\bf (i)}]
$b_1=0$ or $\mbox{\em sgn}(b_1)=\mbox{\em sgn}(a_ma_{m+1}-b_1b_2),$
\item[{\bf (ii)}]
$b_2=0$ or $\mbox{\em sgn}(b_2)=\mbox{\em sgn}(a_ma_{m+1}-b_1b_2).$
\end{description}
\end{cor}

\begin{thm}
\label{czclass}
Let $C\in\R^{n,n}$ be a nonsingular cyclops with eye $m+$ and parameters
$a_i, b_1, b_2$, whose inverse is a Z--matrix. Let
$$ \chi=\{\ k-j\ |\ a_ja_k-b_1b_2<0,\ j\leq m,\ k\geq m+1\ \},$$
$$ \Upsilon=\{\ k-j\ |\ a_ja_k-b_1b_2>0,\ j\leq m,\ k\geq m+1\ \},$$
and  define
\[
x=\left\{ \begin{array}{l}
\min(\chi) \ \ \ \mbox{if $\chi\neq\emptyset$} \\
n+1 \ \ \ \mbox{otherwise,}
\end{array}\right.
\]
\[
y=\left\{ \begin{array}{l}
\min(\Upsilon)\ \ \ \mbox{if $\chi\neq\emptyset$} \\
n+1 \ \ \ \mbox{otherwise,}
\end{array}\right.
\]
$$r=\mbox{number of positive $a_j$ with $j\leq m$},$$
$$t=\mbox{number of positive $a_j$ with $j> m$}.$$
Then the following hold:
\begin{description}
\item[{\bf (i)}] If \mbox{\em det}$C>0$ and $a_m>0$, then $C\in L_n$
(i.e., C is an inverse M--matrix.)
\item[{\bf (ii)}] If \mbox{\em det}$C>0$ and $a_m\leq 0$, then $C\in L_s$,
where $$s=n-1-\max\{m,n-m,n-x+1\}.$$
\item[{\bf (iii)}] If \mbox{\em det}$C<0$, then $C\in L_s$, where
$$s=n-1-\max\{r,t,n-y+1\}.$$
\end{description}
\end{thm}

\pf
Let $B$ be any principal submatrix of $C$, partitioned as in (\ref{cyclops}).
Then $B$ is one of three types:
\begin{enumerate}
\item
$B$ is a principal submatrix of $C_{11}$, in which case,
by (\ref{detd}), det$B$ has
the same sign as the parameter $a_j$ with the largest index contained in $B$.
\item
$B$ is a principal submatrix of $C_{22}$, in which case det$B$ has the same
sign as the parameter $a_j$ with the smallest index contained in $B$.
\item
$B$ is neither a principal submatrix of $C_{11}$ nor of $C_{22}$;
in this case det$B$ has
the same sign as $a_ja_k-b_1b_2$, where $a_j$ has the largest index
less than $m$ contained in $B$, and $a_k$ has the smallest index greater
than $m+1$ contained in $B$.
\end{enumerate}

If \mbox{det}$C>0$, by Corollary \ref{singcyc} we have that
$a_ma_{m-1} - b_1b_2 > 0$. Since $C$ is an inverse Z--matrix,
it follows from our previous results that $b_1 \geq 0, b_2\geq 0$.
Hence $a_m$ and $a_{m+1}$ are nonzero and have the same sign. If
$a_m>0$, it follows that $C$ is nonnegative and hence an inverse M--matrix,
i.e., (i) holds.
If $a_m\leq0$, then $a_{m+1}<0$ and $\mbox{det}C_{11}, \ \mbox{det}C_{22}$
are both negative; thus $C$ has negative principal minors of sizes $m\times m$
and $(n-m)\times (n-m)$.
We need also consider submatrices of the third type; the
largest such submatrix with a negative determinant is of size
$(n-x+1)\times (n-x+1)$. By Theorem \ref{lsz} applied to $A=C^{-1}$,
we have that (ii) holds.

If \mbox{det}$C<0$, then $s$ is determined by the size of
the largest submatrix of $C$ with a positive determinant.
The largest submatrix of $C_{11}$ with a positive determinant is $r\times r$.
The largest submatrix of $C_{22}$ with a positive determinant is $t\times t$.
The largest submatrix of $C$ of the third type is of size
$(n-y+1)\times (n-y+1)$, and (iii) follows.
\eop

\begin{exa}
{\em
The following example illustrates a cyclops and its tridiagonal inverse,
computed by Theorem \ref{invcyc}.
\begin{eqnarray*}
C = \left[ \begin{array}{cccccc}
 4 &    3  &   2  &  -1 &   -1 &   -1 \\
     3  &   3  &   2  &  -1   & -1  &  -1 \\
     2 &    2 &    2  &  -1 &   -1  &  -1 \\
    -4 &   -4  &  -4   &  1   &  1  &   1 \\
    -4  &  -4  &  -4  &   1  &   2  &   2 \\
    -4 &   -4  &  -4  &   1  &   2   &  3
\end{array} \right],
\ \ \
C^{-1} = \left[ \begin{array}{cccccc}
    1 &   -1 &  & & &\\
   -1 &  2 &   -1 &   & &\\
        &   -1  &  0.5 &  -0.5  & & \\
       & &    -2 &        0 &  -1  & \\
      & & & -1 &     2 &   -1 \\
     & & & & -1 &    1
\end{array} \right].
\end{eqnarray*}
Note that, as shown in Theorem \ref{invcyc}, all row sums and column sums
of $C^{-1}$ are zero expect the $3$--rd and the $4$--th. Moreover, using
Theorem
\ref{detcyc}, one easily obtains that det$C=-2$. Applying Theorem
\ref{czclass}, we have $y=3,r=3,t=3$, and thus $C^{-1} \in L_1.$
}
\end{exa}

Next, we will characterize generalized ultrametric matrices whose inverses are
tridiagonal. We begin with the irreducible case. We remind the reader that
$C\in\R^{n,n}$ is called {\em irreducible} if its directed graph, $\Gamma(C)$,
is strongly connected (see e.g., \cite{BP}). Also recall that $C$ is called
{\em totally nonnegative} if all its minors are nonnegative.

\begin{thm}
\label{igum}
Let $C\in\R^{n,n}$ be a nonsingular matrix. Then the following are equivalent:
\begin{description}
\item[{\bf (i)}]
There is $m \in \l n\ra$
such that $C^{-1}$ is an irreducible row and column diagonally dominant
tridiagonal M--matrix whose row and columns sums are all zero,
except at least one of the $m$--th or $(m+1)$--st.
\item[{\bf (ii)}] There is $m \in \l n\ra$ and $a_i,b_1,b_2\in\R$ such
that $C$ is a cyclops with eye $m+$ and parameters $a_i,b_1,b_2$ satisfying
$$
\min\{a_m,a_{m+1}\}\geq \max\{b_1,b_2\} \ \ \mbox{and} \ \ \
\min\{b_1,b_2\}>0.
$$
\item[{\bf (iii)}] $C$ is an irreducible generalized ultrametric matrix
whose inverse is tridiagonal.
\item[{\bf (iv)}]
$C$ is a totally nonnegative irreducible generalized ultrametric matrix.
\end{description}
\end{thm}

\pf
Let $C=[c_{ij}]$ and $A=C^{-1}$.

{\bf (i) implies (ii):}
By \cite[Theorem 3.2]{MNST1} applied to $A$ and $A^T$,
$c_{ii}=c_{ij}=c_{ji}$ for all $1\leq j<i\leq m$ and $c_{ii}=c_{ik}=c_{ki}$
for all
$m+1\leq i< k\leq n$. Since $A$ is a tridiagonal M--matrix,
by \cite[Theorem 4.1]{MNST2},
$c_{mk}= \frac{c_{m,m+1}c_{m+1,k}}{c_{m+1,m+1}}=c_{m,m+1}$ and
$c_{jk}=\frac{c_{jm}c_{mk}}{c_{mm}}=c_{mk}$ for all $j\leq m,$ and $k\geq
m+1$. Similarly, $c_{jk}=c_{m+1,k}=c_{m+1,m}$ for all $j\geq m+1, k\leq m$.
Thus $C$ is a cyclops with eye $m+$ and parameters
$a_i=c_{ii},b_1=c_{m,m+1},b_2=c_{m+1,m}$.
Since $A$ is a nonsingular M--matrix, all principal minors of $C$ are
positive and thus $a_{m}a_{m+1}-b_1b_2>0$. Since $A$
is irreducible, the entries of $C$ are all positive. Lastly, the inequality
$min\{a_m,a_{m+1} \}\geq max\{b_1,b_2\}$ follows from the facts that $C$ is an
entrywise positive cyclops and that $C^{-1}$ is a row and column diagonally
dominant matrix (and thus each diagonal entry of $C$ is greater than
or equal to the other entries in the corresponding row and column, see
\cite[Theorem (3,5)]{FP}).

{\bf (ii) implies (iii):}
Notice first that the conditions on the parameters of the cyclops imply
that $C$ is an generalized ultrametric matrix. By Corollary \ref{zm},
$C$ is invertible, and by \cite[Theorem 4.1]{MNST2}, $C^{-1}$ is tridiagonal.
%Follows from Corollary \ref{zm} and the fact that a
%Z--matrix whose inverse is a nonnegative matrix must be an M--matrix
%(see \cite{BP}).

{\bf (iii) implies (i):}
We only need to show that the row and column sums are as claimed.

\begin{description}
\item[Claim I:] For all $j<i<k$, either $c_{ii}=c_{ik}=c_{ki}$ or
$c_{ii}=c_{ij}=c_{ji}$.

Proof of Claim I:
>From \cite[Theorem 4.1]{MNST2}, we see that for all $j<i<k$,
$c_{jk}=\frac{c_{ji}c_{ik}}{c_{ii}}$ and $c_{kj}=\frac{c_{ki}c_{ij}}{c_{ii}}.$
Consider the {\em triangle} on vertices $i,j,k$ (see \cite[Definition
2.3]{MNST1}).
If
$i$ is the {\em preferred vertex}, then $c_{ji}=c_{ki}$ and
$c_{ij}=c_{ik}$, and hence
$c_{jk}=c_{kj}$ and $c_{jk}\geq max \{c_{ji},c_{ik}\}$. But since $C$ is a
generalized ultrametric matrix,
$c_{ii} \geq max \{c_{ji},c_{ik}\}$, thus
$c_{ii}=c_{ij}=c_{ji}=c_{jk}=c_{kj}=c_{ik}=c_{ki}$. If $j$ is preferred, it
follows
that $c_{jk}=c_{ji}$ and $c_{kj}=c_{ij}$ and thus $c_{ii}=c_{ik}=c_{ki}$. If
$k$ is preferred then $c_{ii}=c_{ij}=c_{ji}$. This establishes Claim I.

\item[Claim II:] If $p$ is the first nonzero row sum of A, then all other
row sums
are zero, except possibly the $(p+1)$--st.

Proof of Claim II: For all $j<p+1,$ \cite[Theorem 3.2]{MNST1}
implies $c_{p+1,p+1}\neq c_{j,p+1}$.
Suppose there is $q>p+1$ such that $q$--th row sum is nonzero.
Then for all $k>p+1,$ \cite[Theorem 3.2]{MNST1} implies $c_{p+1,p+1}\neq
c_{k,p+1}$.
If we apply Claim I with $i=p+1$ we have a contradiction that establishes
Claim II.

\item[Claim III:] If the $q$--th row sum of $A$ is nonzero, then the column
sums
of columns $1,2,\ldots, q-2$ are zero.

Proof of Claim III:
By \cite[Theorem 3.2]{MNST1},
for all $i<q$ and $k\geq q$,
$c_{ii}\neq c_{ki}$ and hence by Claim I, for all $j<i$,
$c_{ii}=c_{ij}=c_{ji}$. If we now apply \cite[Theorem 3.2]{MNST1} to $A^T$,
we see that columns
$1,2,\ldots, q-1$ must have zero column sums. This establishes Claim III.
\end{description}
Suppose now that the $p$--th row sum of $A$ is the first nonzero row sum.
Let $P$ be the permutation matrix which reverses the order of the indices
$1,2,\ldots, n$.

If the $(p+1)$--st row sum of $A$ is
also nonzero, then by applying Claim III to $A$ and $PAP^T$, with $q=p$ and
with
$q=p+1$, we see that the only possible nonzero column sums are the $p$--th and
the $(p+1)$--st. Taking $m=p$, the implication is proven.

If the $(p+1)$--st row sum is zero, then
by applying Claim III to $A$ and $PAP^T$, with $q=p$, we see that the
only possible
nonzero column sums are the $(p-1)$--st, the $p$--th and the $(p+1)$--st.
By Claim II
applied to $A^T$, either the $(p-1)$--st or the $(p+1)$--st sum is zero.
By choosing
$m$ appropriately to be either $p-1$ or $p$, the implication is proven.

{\bf (iii) if and only if (iv):}
Follows from the results in \cite{MNST1} or \cite{NV}, and in
\cite{L80}. \eop

\begin{cor}
A matrix $A$ is of type $\D$ with parameter $a_1>0$ if and
 only if $A^{-1}$ is a tridiagonal M--matrix with the only nonzero row and
column
sums being the first. \end{cor}

We say that $C$ is a {\em G--cyclops} if it is
nonsingular and satisfies any of the equivalent conditions of Theorem
\ref{igum}. We
also refer to a matrix all of whose entries are equal as a {\em flat} matrix.

\begin{thm}
\label{rgum}
Let $C$ be a nonsingular matrix that is reducible
but not completely reducible. Then the following are equivalent:
\begin{description}
\item[{\bf (i)}] Either $C$ or $C^T$ is of the form
\begin{displaymath}
B := \pmatrix{B_{11} & B_{12} & B_{13} & \ldots & B_{1m} \cr
0 & B_{22} & B_{23} & \ldots & B_{2m} \cr
\vdots & \ddots & \ddots & \vdots & \vdots \cr
\vdots & &\ddots & B_{m-1,m-1} & B_{m-1,m}\cr
0 & \ldots & \ldots & 0 & B_{mm} \cr},
\end{displaymath}
where $B_{11}$ is a G--cyclops whose last
column is of constant value given by $f_{11}$; $B_{mm}$ is a G--cyclops
whose first row is of constant value given by $f_{mm}$;
each $B_{st},$ $1\leq s < t \leq m$, is a
flat matrix whose fixed value is given by $f_{st}$; each
$B_{ss},$ $2\leq s \leq m-1$, is either a positive number $f_{ss}$ or
an entrywise positive matrix of the form
\begin{displaymath}
B_{ss} \ = \
\pmatrix{f_{ss} & f_{ss} \cr
g_{ss} & f_{ss} \cr},
\end{displaymath}
with $g_{ss} < f_{ss}$; and for
some $2\leq q\leq m-1$, $$f_{11}\geq f_{22}\geq \ldots\geq f_{qq}>0,$$
$$0<f_{q+1,q+1}\leq f_{q+2,q+2}\leq\ldots\leq f_{mm},$$
$$f_{1t}=f_{2t}=\ldots=f_{t-1,t}=f_{tt} \ \ \mbox{for} \ \ 1< t\leq q,$$
$$f_{ss}=f_{s,s+1}=\ldots=f_{s,m-1}=f_{sm} \ \ \mbox{for} \ \ q< s<m,$$
$$0<f_{st}=f_{q,q+1}\leq \min\{f_{qq},f_{q+1,q+1}\} \ \ \mbox{for} \ \
1\leq s\leq
q<t\leq m.$$
\item[{\bf (ii)}]
$C$ is a generalized ultrametric matrix \ whose
inverse is tridiagonal. \item[{\bf (iii)}]
$C$ is a totally nonnegative generalized ultrametric matrix.
\end{description}
\end{thm}

\pf

{\bf (i) implies (ii):}
It is easy to see that the conditions on the
$f_{ij}$ guarantee that $C$ is a generalized ultrametric matrix and that it
satisfies \cite[Theorem 4.1 (ii)]{MNST2}.

{\bf (ii) implies (i):}
Let $C=[c_{ij}]$ and $A=C^{-1}=[a_{ij}]$.
Since $A$ is a tridiagonal M--matrix, it must satisfy
\cite[Theorem 4.1 (ii)]{MNST2}.

We begin by showing that $C$ or $C^T$ must
be block upper triangular with no zero entries in or above the diagonal
blocks. Suppose that
$A$ has a zero entry on the superdiagonal and a zero entry on the subdiagonal.
For simplicity, we will assume that
$a_{j,j+1}=0$ and $k\geq j$ is the smallest integer
such that $a_{k+1,k}=0$ (otherwise take $A=(C^T)^{-1}).$ If $k=j$ then $A$ is
completely reducible contradicting our hypothesis, hence we will assume
that $j< k$.
Since $j$ does not access $j+1$ in $\Gamma(C)$, by \cite[Lemma 2.2]{S},
$c_{j,j+1}=0.$ Similarly $c_{k+1,k}=0$. By \cite[Theorem 4.1 (ii)]{MNST2},
$c_{ji}=\frac{c_{j,j+1}c_{j+1,i}}{c_{j+1,j+1}}=0$ for all $i\geq j+1$, and
$c_{k+1,l}=\frac{c_{k+1,k}c_{kl}}{c_{kk}}=0$ for all $l\leq k$.
Consider the triangle
(see \cite[Definition 2.3]{MNST2}) on $j,k,k+1$. Either $c_{kj}=0$
or $c_{k,k+1}=0$.
If $c_{k,k+1}=0$ then by \cite[Lemma 2.2]{S} $a_{k,k+1}=0$ and hence $A$ is
completely reducible, contradicting our hypothesis. If $c_{kj}=0$, then by
\cite[Lemma 2.2]{S}, $k$ does not have access to $j$ in $\Gamma(A)$ and
hence there
must be an $i$ with $j\leq i < k$ such that $a_{i+1,i}=0$, contradicting the
minimality of $k$. So either the superdiagonal of $C^{-1}$ or the subdiagonal
of $C^{-1}$ contains only nonzero entries. It follows that either
$C$ or $C^T$ is block
upper triangular, as represented by $B$.
It remains to show that the blocks of $B=[b_{ij}]$ are as claimed.
Without loss of generality, assume that all
the entries on the superdiagonal of $A$ are nonzero.

Since both $B$ and $A$ are block upper triangular,
$B_{ss}=(A_{ss})^{-1}$. Since $A_{ss}$ is irreducible, $B_{ss}$
must be a G--cyclops
by Theorem \ref{igum}. Since the superdiagonal entries of $A$ are nonzero, by
the results in \cite{S} each $B_{st}$ is an entrywise positive matrix
for all $s\leq t$.

Notice that if $i<j$ and $i$ and $j$ are in
different blocks of $B$ then $b_{ji}=0.$ This fact will be used without
further remark whenever triangles are considered in the remainder of this
proof. We will also write $j\in s$ to mean that $b_{jj}$ is in the block
$B_{ss}$.

Let $j\in 1, \ k\in 1$ and $l\notin 1$ with $j\leq k$. By
considering the triangle on
$j,k,l$, we see that $b_{jl}=b_{kl}.$ But
$b_{jl}=\frac{b_{jk}b_{kl}}{b_{kk}}$, hence $b_{jk}=b_{kk}$ and $B_{11}$
is as claimed. A similar argument shows that $B_{mm}$ is as claimed.

Consider $B_{rr},\ B_{ss},$ and $B_{tt}$
with $r<s<t$. Let $i\in r,\ j\in s, \ k\in s,\ l\in t$. Consider the triangle
on
$i,j,l$. Then $b_{il}=min\{ b_{ij}, b_{jl}\}.$ By \cite[Theorem 4.1]{MNST2},
$b_{il}=\frac{b_{ij}b_{jl}}{b_{jj}}$ and hence $b_{jj}=max\{b_{ij},b_{jl}\}$.
Similarly, $b_{kk}=max\{b_{ik},b_{kl}\}$. From the triangle on $i,j,k$
we see that
$b_{ij}=b_{ik}\leq max \{ b_{jk},b_{kj}\}$. From the triangle on $j,k,l$
we have that
$b_{jl}=b_{kl}\leq max \{ b_{jk},b_{kj}\}$. But then
$b_{jj}=max\{b_{ij},b_{jl}\} \leq max\{ b_{jk},b_{kj}\}\leq b_{jj}$.
Hence equality must hold throughout. Using the corresponding inequalities
for $k$, we can conclude that
$b_{jj}=\max\{b_{jk},b_{kj}\}=b_{kk}$. If $j<k$, then $b_{ik}=\frac{
b_{ij}b_{jk}}{b_{jj}}$ implies that $f_{ss}=b_{jk}=b_{jj}\geq b_{kj}=g_{ss}$.
Since $B_{ss}=(A_{ss})^{-1}$, $B_{ss}$ must be a nonsingular G--cyclops and
hence can only be as claimed for $2\leq s\leq m-1$.

Let now $h=min\{ j \ | \ b_{jj}=b_{jk}$ for all $k>j \}$
($h$ is well defined since $b_{nn}=b_{nk}$ for $k>j$).
Consider the $r$--th diagonal
block so that $h\in r$. By the equalities in
the above paragraph, if
$h-1\in r$ then $b_{h-1,h-1} =b_{h-1,h}=b_{hh}=b_{hk}=b_{h-1,k}$,
contradicting the
minimality of $h$. Hence $h-1\notin r$. Set $q=r-1$.
If $r<m$, then by the choice of
$h$ and the triangle on $h,h+1,k,$ $b_{hh}=b_{h,h+1}=b_{hk}=min\{ b_{h,h+1},
b_{h+1,k}\}$ for all $k> h+1,$ which implies that
$b_{h,h+1}\leq b_{h+1,k}.$ Thus
$b_{h+1,h+1}=max\{b_{h,h+1},b_{h+1,k}\}=b_{h+1,k}$ for all $k> h+1.$
We can now repeat this argument for $h+2,h+3,\ldots,$
up to largest index in the $(m-1)$--st
diagonal block to conclude that for all $q<s<t$ with $j\in s$ and $k\in t$,
$f_{ss}=b_{jj}=b_{jk}=f_{st}$ and thus the $B_{st}$ are as claimed.
For any $j<h$, with
$j\notin 1$, by the choice of $h$ and the inequalities in the above
paragraph, we see that
$b_{jj}=b_{ij}$ for all $i\leq j$. Hence $B_{st}$ must be as claimed for
all $s<t\leq q$. For $s\leq q < t$, let $j\in s$ and $k\in t$. Then
$b_{jk}=\frac{b_{jh}b_{hk}}{b_{hh}}=b_{jh}$ and
$b_{jh}=\frac{b_{j,h-1}b_{h-1,h}}{b_{h-1h-1}}$.
By considering the triangle on
$j,h-1,h$, we have that $b_{jh}=min\{ b_{j,h-1},b_{h-1,h}\}=min\{
b_{h-1,h-1},b_{h-1,h}\} =b_{h-1,h}\leq min\{ b_{h-1,h-1},b_{hh}\}.$
Thus $f_{st}=b_{jk}=b_{h-1,h}=f_{q,q+1}.$

{\bf (ii) if and only if (iii):}
Follows from the results in \cite{MNST1}
or \cite{NV}, and in \cite{L80}. \eop

\begin{exa}
{\em The following matrix illustrates a matrix that satisfies the
conditions of Theorem \ref{rgum} (and hence it is a reducible,
totally nonnegative, generalized ultrametric matrix whose inverse is
tridiagonal).

\begin{eqnarray*}
C = \left[ \begin{array}{rrrrrrrrrrrr}
 12 & 11 & 10 & 9 & 9 & 7 & 5 & 5 &    5 &    5 &    5 &    5 \\
 11 &   11 &   10  &   9 &    9 &    7 &    5 &    5 &    5 &    5 & 5 & 5\\
 10 &   10 &   10  &   9 &    9 &    7 &    5 &    5 &    5 &    5 & 5 & 5\\
  0 &    0 &    0  &   9 &    9 &    7 &    5 &    5 &    5 &    5 & 5 & 5\\
  0 &    0 &    0  &   8 &    9 &    7 &    5 &    5 &    5 &    5 & 5 & 5 \\
  0 &    0 &    0  &   0 &    0 &    7 &    5 &    5 &    5 &    5 & 5 & 5 \\
  0 &    0 &    0  &   0 &    0 &    0 &    6 &    6 &    6 &    6 & 6 & 6 \\
  0 &    0 &    0  &   0 &    0 &    0 &    0 &    7 &    7 &    7 & 7 & 7 \\
  0 &    0 &    0  &   0 &    0 &    0 &    0 &    6 &    7 &    7 & 7 & 7\\
  0 &    0 &    0  &   0 &    0 &    0 &    0 &    0 &    0 &    8 & 8 & 8\\
  0 &    0 &    0  &   0 &    0 &    0 &    0 &    0 &    0 &    8 & 9 & 9 \\
  0 &    0 &    0  &   0 &    0 &    0 &    0 &    0 &    0 &    8 & 9 & 10
\end{array} \right].
\end{eqnarray*}
}
\end{exa}
Finally, Theorems \ref{igum} and \ref{rgum} yield the following result.
\begin{thm}
\label{ggum}
Let $C\in\R^{n,n}$ be a nonsingular matrix.
Then the following are equivalent:
\begin{description}
\item[{\bf (i)}]
$C$ is the direct sum of matrices of the forms given in
Theorem \ref{igum}(ii) and Theorem \ref{rgum}(i).
\item[{\bf (ii)}]
$C$ is a generalized
ultrametric matrix whose inverse is tridiagonal.
\item[{\bf (iii)}]
$C$ is a totally nonnegative generalized ultrametric matrix.
\end{description}
\end{thm}

%****************************************************************
%	BIBLIOGRAPHY
%****************************************************************

\end{document}